\documentclass[reqno]{amsart}

\usepackage{amsmath,amssymb,amsthm,graphicx,a4wide,enumerate}
\usepackage[small,bf]{caption}
\usepackage{hyperref}
\theoremstyle{plain}

\theoremstyle{definition}

\theoremstyle{remark}

\newcommand{\prn}[1]{\left(#1\right)}

\newcommand{\brk}[1]{\left[#1\right]}
\newcommand{\brc}[1]{\left\{#1\right\}}

\newcommand{\pd}[2]{\frac{\partial#1}{\partial#2}}
\newcommand{\ud}[1]{\,\text{d}#1}

\begin{document}
\parskip1ex

\setlength{\oddsidemargin}{0.5\paperwidth}
\addtolength{\oddsidemargin}{-0.5\textwidth}
\addtolength{\oddsidemargin}{-1in}
\setlength{\evensidemargin}{\oddsidemargin}
\setlength{\textheight}{8.6in}
\setlength{\topmargin}{-0.1in}
\setlength{\headsep}{.2in}
\setlength{\footskip}{.3in}


\title[A particle method for traffic flow simulations on highway networks]{A characteristic particle method for traffic flow simulations on highway networks}

\author{Yossi Farjoun}
\address[Yossi Farjoun]
{G.\ Mill\'an Institute of Fluid Dynamics \\ Nanoscience and Industrial Mathematics \\
University Carlos III de Madrid \\ Av.\ Universidad 30, 28911 Legan\'es \\ Spain}
\email{yfarjoun@ing.uc3m.es}
\author{Benjamin Seibold}
\address[Benjamin Seibold]
{Department of Mathematics \\ Temple University \\
1805 North Broad Street \\ Phila\-delphia, PA 19122}
\email{seibold@temple.edu}
\urladdr{http://www.math.temple.edu/\~{}seibold}
\date{\today}

\maketitle

\begin{abstract}
A characteristic particle method for the simulation of first order macroscopic traffic models on road networks is presented. The approach is based on the method \emph{particleclaw}, which solves scalar one dimensional hyperbolic conservation laws exactly, except for a small error right around shocks. The method is generalized to nonlinear network flows, where particle approximations on the edges are suitably coupled together at the network nodes. It is demonstrated in numerical examples that the resulting particle method can approximate traffic jams accurately, while only devoting a few degrees of freedom to each edge of the network.
\end{abstract}


\section{Introduction}
Macroscopic traffic models describe the evolution of the vehicular traffic density on a road by partial differential equations. A large network of highways, which includes ramps and intersections, can be modeled as a directed graph, whose edges are road segments that join and bifurcate at network nodes. The traffic density on each edge evolves according to some macroscopic traffic model, and at nodes, specific coupling conditions are specified. Since realistic road networks can easily consist of tens on thousands of edges, efficient numerical methods are crucial for the simulation, forecasting, now-casting, or optimization of traffic flow on networks.

In this paper, we present a characteristic particle method for ``first order'' traffic models on road networks. The approach is based on the method \emph{particleclaw} \cite{FarjounSeibold:FarjounSeibold2009, FarjounSeibold:Particleclaw}, which solves scalar one dimensional hyperbolic conservation laws exactly, except for the immediate vicinity of shocks, where a small approximation error is observed. Particleclaw is used to evolve the numerical solution on each edge, and a methodology is presented for the coupling of the particle approximations on individual edges together through the network nodes. In recent years \cite{FarjounSeibold:FarjounSeiboldMeshfree2008, FarjounSeibold:FarjounSeibold2010, FarjounSeibold:FarjounSeiboldMeshfree2011}, particleclaw has been demonstrated to possess certain structural advantages over approaches that are based on a fixed grid. To name a few examples:
\begin{itemize}
\item It possesses no numerical viscosity, and thus preserves a reasonable accuracy even with very few degrees of freedom. This is in contrast to Godunov's method \cite{FarjounSeibold:Godunov1959} and other low order schemes.
\item It is optimally local, in the sense that particles move independently of each other, and communication between neighboring particles occurs only at shocks. In contrast, high order finite difference approaches \cite{FarjounSeibold:LaxWendroff1960}, finite volume methods \cite{FarjounSeibold:VanLeer1974}, or ENO \cite{FarjounSeibold:HartenEngquistOsherChakravarthy1987}/WENO \cite{FarjounSeibold:LiuOsherChan1994} schemes use wide stencils to achieve high order, which poses a particular challenge right at network nodes.
\item It is total variation diminishing, yet it is second order accurate even in the presence of shocks \cite{FarjounSeibold:FarjounSeibold2009}. In contrast, many limiters in fixed grid approaches need to reduce the order of convergence near shocks in order to avoid overshoots.
\item The approximate solution is represented as a piecewise similarity solution of the underlying equation, that is continuous except for right at actual shocks. In contrast, the reconstruction in finite volume \cite{FarjounSeibold:VanLeer1974} and Discontinuous Galerkin \cite{FarjounSeibold:ReedHill1973, FarjounSeibold:CockburnShu1988} methods tends to possess spurious discontinuities at every boundary between cells.
\item It is adaptive by construction. Particles adapt to the shape of the solution, and the approach can be generalized to incorporate stiff source terms \cite{FarjounSeibold:FarjounSeiboldMeshfree2011}. This is in contrast to fixed grid methods, which would need to use mesh refinement techniques \cite{FarjounSeibold:BergerOliger1984} for adaptivity.
\end{itemize}
These advantages render \emph{particleclaw} as an interesting candidate for the simulation of traffic flow on networks. The need for the design of specialized numerical schemes for traffic models on networks has been pointed out for instance by Bretti, Natalini, and Piccoli \cite{FarjounSeibold:BrettiNataliniPiccoli2006b}, who design an efficient numerical approach for a very specific class of flux functions, for which it suffices to track the location of the transition between free and congested traffic flow. The idea of tracking features is also shared by front tracking methods \cite{FarjounSeibold:HoldenHoldenHeghKrohn1988}, which approximate the whole solution by finitely many moving discontinuities. In fact, in analogy to front tracking, \emph{particleclaw} can also be interpreted as a \emph{rarefaction tracking} \cite{FarjounSeibold:FarjounSeibold2010}.

This paper is organized as follows. In Sect.~\ref{FarjounSeibold:sec:traffic_models_networks}, the class of traffic models under consideration is outlined, and the coupling conditions for road networks are given. The characteristic particle method \emph{particleclaw} is presented in Sect.~\ref{FarjounSeibold:sec:particleclaw}, and in Sect.~\ref{FarjounSeibold:sec:particleclaw_networks}, we demonstrate how the approach can be generalized to nonlinear flows on networks. Numerical results are shown in Sect.~\ref{FarjounSeibold:sec:numerical_results}, and in Sect.~\ref{FarjounSeibold:sec:conclusions_outlook}, we present conclusions and give an outlook.

\section{Traffic Models on Highway Networks}
\label{FarjounSeibold:sec:traffic_models_networks}

\subsection{Macroscopic Traffic Models}
Macroscopic traffic models treat the traffic as a continuum, and use partial differential equations to describe the temporal evolution of the vehicle density $u(x,t)$,\footnote{Even though densities are commonly denoted by $\rho$, here we use $u$, in order to express the fact that the numerical approach applies to more general network flows.} where $x$ is the position on the road (the flow is averaged over all lanes that go in one direction), and $t$ is time. If vehicles move with the velocity $v(x,t)$, then the conservation of vehicles (in the absence of ramps) is described by the continuity equation $u_t+(uv)_x = 0$. The assumption of a direct functional relationship between the velocity and the density, $v = v(u)$, yields the classical Lighthill-Whitham-Richards (LWR) model \cite{FarjounSeibold:LighthillWhitham1955, FarjounSeibold:Richards1956}
\begin{equation}
\label{FarjounSeibold:eq:lighthill_whitham}
u_t+\prn{f(u)}_x = 0\;,
\end{equation}
which is a scalar hyperbolic conservation law with a flux function $f = uv$ that equals the traffic flow rate. Due to the direct density-velocity relationship, the LWR model does not model vehicle accelerations, and thus is not able to describe non-equilibrium features such as phantom traffic jams or self-sustained traffic waves (``jamitons'') \cite{FarjounSeibold:FlynnKasimovNaveRosalesSeibold2009}. To overcome this limitation, numerous types of ``second order'' models have been developed, such as the Payne-Whitham model \cite{FarjounSeibold:Payne1979}, the Zhang-Aw-Rascle model \cite{FarjounSeibold:Zhang1998, FarjounSeibold:AwRascle2000}, phase transition models \cite{FarjounSeibold:Colombo2003}, and many more. These introduce a velocity (or velocity-like) variable as an independent unknown into the equations, resulting in a system of balance laws. In this paper, we restrict our study to effects that can be captured reasonably well by the scalar LWR model, i.e.~we are interested in the large scale, nonlinear equilibrium behavior of traffic flow.

As it is common in traffic models \cite{FarjounSeibold:Helbing2001}, we assume that the LWR flux function satisfies the following conditions: (i) no flow for empty road and bumper-to-bumper densities, i.e.~$f(0) = 0 = f(u^\text{m})$, for some maximum density $u^\text{m}$; and (ii) concavity, i.e.~$f''(u)<0\;\forall\,u\in[0,u^\text{m}]$. As a consequence of these assumptions, there is a unique maximum flow rate $f^* = f(u^*)$ that occurs for an optimal density $u^*$. Two parameters of $f$ are fundamentally important: the slope $v^\text{m} = f'(0)$ is the velocity of vehicles when alone on the road (i.e.~approximately the speed limit); and $u^\text{m}$ is the number of lanes divided by the average vehicle length plus a safety distance. Other than for these two values, the precise functional shape of $f$ can depend on many factors, such as the road geometry, the drivers' state of mind, etc. Frequently, one assumes a simple parabolic profile $f(u) = v^\text{m}u\,(1-\frac{u}{u^\text{m}})$. This particular shape was inspired by measurements of Greenshields \cite{FarjounSeibold:Greenshields1935}, and even though it does not fit well with contemporary measurements \cite{FarjounSeibold:Helbing2001}, it is commonly used due to its simplicity. We shall do so here as well, as it simplifies the presentation of the numerical approach in Sect.~\ref{FarjounSeibold:sec:particleclaw_networks}.

\subsection{Traffic Networks}
\label{FarjounSeibold:subsec:traffic_networks}
A traffic network is a directed graph of \emph{network edges} (road segments) that join and bifurcate at \emph{network nodes}. On each edge, the traffic evolves according to the LWR model \eqref{FarjounSeibold:eq:lighthill_whitham} with a flux function that is specific to this edge (here, we assume that on each edge the flux function does not explicitly depend on space or time). In order to have a feasible model for traffic flow on road networks, one must formulate coupling conditions that guarantee the existence of a unique entropy solution, given suitable initial and boundary conditions. This problem was first addressed by Holden and Risebro \cite{FarjounSeibold:HoldenRisebro1995}, and subsequently generalized by various other authors (see below). In the following, we briefly outline the key ideas of coupling conditions at network nodes.

At a network node, one has the following setup. Let the node have $n$ ingoing edges (roads), numbered $i = 1,2,\dots,n$, each of which carries a vehicle density $u_i(x_i,t)$, where $x_i\in [0,L_i]$, and a flux function $f_i(u)$. Similarly, let the node have $m$ outgoing edges, numbered $j = n+1,\dots,n+m$, each carrying a density $u_j(x_j,t)$, where $x_j\in [0,L_j]$, and a flux function $f_j(u)$. Ingoing edges end at the node at $x_i = L_i$, and outgoing edges start at the node at $x_j = 0$. The conservation of vehicles requires that the total inflow into the node equals the total outflow out of the node, i.e.
\begin{equation}
\label{FarjounSeibold:eq:node_conservation}
\sum_{i=1}^n f_i(u_i(L_i,t)) = \sum_{j=n+1}^{n+m} f_j(u_j(0,t))\;\forall\,t\;.
\end{equation}
In order to obtain further rules for the temporal evolution of the
solution at nodes, a generalized Riemann problem is formulated, as
follows. Let all edges be extended away from the node to $\pm\infty$,
and consider initial conditions that are constant on each edge,
i.e.~$u_i(x,0) = u_i$, for $x\in\brk{-\infty,L_i}$, and $u_j(x,0) = u_j$
for $x\in\brk{0,\infty}$. The question is then: what are the new states $\hat{u}_i = u_i(L_i,t)$ and $\hat{u}_j = u_j(0,t)$ at the node for $t>0$? These new states define the solution at all times, since the problem admits a self-similar solution $u_i(x,t) = w_i(\frac{L_i-x}{t})$ and $u_j(x,t) = w_j(\frac{x}{t})$. Let the flux values of the new states be denoted $\gamma_i = f_i(\hat{u}_i)\;\forall\,i\in\{1,\dots,n+m\}$. Clearly, these fluxes must satisfy the conservation condition \eqref{FarjounSeibold:eq:node_conservation}, i.e.
\begin{equation}
\label{FarjounSeibold:eq:node_conservation_gamma}
\sum_{i=1}^n \gamma_i = \sum_{j=n+1}^{n+m} \gamma_j\;.
\end{equation}
Moreover, the new states $\hat{u}_i\;\forall\,i\in\{1,\dots,n+m\}$ must generate a solution for which information on each edge is either stationary or moves away from the node (information cannot move ``into'' the node). As shown in \cite{FarjounSeibold:HoldenRisebro1995}, this information direction requirement implies that the new fluxes satisfy the inequality constraints
\begin{align}
\label{FarjounSeibold:eq:condition_inequality_ingoing}
\gamma_i\in\Omega_i\;
&\text{, where~~}\Omega_i = [0,f_i(\min\{u_i,u_i^*\})]\;
\forall\,i\in\{1,\dots,n\}\;\text{, and} \\
\label{FarjounSeibold:eq:condition_inequality_outgoing}
\gamma_j\in\Omega_j\;
&\text{, where~~}\Omega_j = [0,f_j(\max\{u_j,u_j^*\})]\;
\forall\,j\in\{n+1,\dots,n+m\}\;.
\end{align}
In words: on ingoing edges, the new flux must be less than or equal to the \emph{demand flux}, which is given by the old flux if $u_i\le u_i^*$, and the maximum flux $f_i^*$ if $u_i>u_i^*$. Likewise, on outgoing edges, the new flux must be less than or equal to the \emph{supply flux}, which is given by the old flux if $u_j\ge u_j^*$, and the maximum flux $f_j^*$ if $u_j<u_j^*$.

Conditions \eqref{FarjounSeibold:eq:node_conservation_gamma}, \eqref{FarjounSeibold:eq:condition_inequality_ingoing}, and \eqref{FarjounSeibold:eq:condition_inequality_outgoing} allow infinitely many possible new fluxes. Further conditions must be provided to single out a unique solution. A reasonable set of additional conditions is given by the drivers' desired destinations matrix, defined as follows. Let $\vec{\gamma}_\text{in}$ denote the vector of ingoing fluxes, and $\vec{\gamma}_\text{out}$ be the vector of outgoing fluxes. The fact that very specific percentages of drivers that enter the node through edge $i$ exit the node on edge $j$ (different drivers have different destinations) can be encoded via the linear system
\begin{equation}
\label{FarjounSeibold:eq:driver_destinations}
A \cdot \vec{\gamma}_\text{in} = \vec{\gamma}_\text{out}\;,
\end{equation}
where
\begin{equation*}
A = \begin{pmatrix}
a_{1,n+1} & \hdots & a_{n,n+1} \\
\vdots & \ddots & \vdots \\
a_{1,n+m} & \hdots & a_{n,n+m}
\end{pmatrix},
\quad
\vec{\gamma}_\text{in} =
\begin{pmatrix} \gamma_1 \\ \vdots \\ \gamma_n \end{pmatrix},
\quad
\vec{\gamma}_\text{out} =
\begin{pmatrix} \gamma_{n+1} \\ \vdots \\ \gamma_{n+m} \end{pmatrix}.
\end{equation*}
The matrix $A$ is column-stochastic, i.e.~all $a_{i,j}\in [0,1]$ and $\sum_{j=n+1}^{n+m}a_{i,j} = 1$. Thus, the flux balance condition \eqref{FarjounSeibold:eq:node_conservation_gamma} is automatically guaranteed by relation \eqref{FarjounSeibold:eq:driver_destinations}, since $\vec{e}^T \cdot \vec{\gamma}_\text{out} = \vec{e}^T \cdot A \cdot \vec{\gamma}_\text{in} = \vec{e}^T \cdot \vec{\gamma}_\text{in}$, where $\vec{e} = (1,\dots,1)^T$. Condition \eqref{FarjounSeibold:eq:driver_destinations} together with the constraints \eqref{FarjounSeibold:eq:condition_inequality_ingoing} and \eqref{FarjounSeibold:eq:condition_inequality_outgoing} yields that $\vec{\gamma}_\text{in}$ must lie in the \emph{feasibility domain}
\begin{equation}
\label{FarjounSeibold:eq:feasibility_domain}
\Omega = \brc{\vec{\gamma}\in\Omega_1\times\dots\times\Omega_n \ | \
  A\cdot\vec{\gamma}\in\Omega_{n+1}\times\dots\times\Omega_{n+m}}
\subset \mathbb{R}^n\;. 
\end{equation}
The selection of a specific $\vec{\gamma}_\text{in}\in\Omega$ requires further modeling arguments. Possible criteria are, for instance, entropy arguments \cite{FarjounSeibold:HoldenRisebro1995}, the modeling of the intersection geometry \cite{FarjounSeibold:HertyKlar2003}, or simply the assumption that drivers behave such that the throughput through the node is maximized \cite{FarjounSeibold:CocliteGaravelloPiccoli2005, FarjounSeibold:BrettiNataliniPiccoli2006a}, i.e.~one solves the linear program
\begin{equation}
\label{FarjounSeibold:eq:node_optimization}
\max\, \vec{e}^T \cdot \vec{\gamma}_\text{in}\;\;
\text{s.t.}\;\vec{\gamma}_\text{in}\in\Omega\;.
\end{equation}
In the examples presented in Sect.~\ref{FarjounSeibold:sec:numerical_results}, we shall follow the latter option, even though the other alternatives are possible as well. The modeling has to be augmented by one small but important detail. It is possible that \eqref{FarjounSeibold:eq:node_optimization} does not possess a unique solution, namely if the extremal boundary of $\Omega$ is perpendicular to $\vec{e}$. In this case, one can introduce the additional constraint that $\vec{\gamma}_\text{in}\in\vec{c}\mathbb{R}$, where $\vec{c}\in\mathbb{R}^n$ is a given constant that models the merging behavior at the node.

The definition of the generalized Riemann problem is finalized by selecting new states as follows. On ingoing roads, choose $\hat{u}_i = u_i$ if $\gamma_i = f_i(u_i)$, otherwise choose $\hat{u}_i\ge u_i^*$, s.t.~$f_i(\hat{u}_i) = \gamma_i$. Similarly, on outgoing roads, choose $\hat{u}_j = u_j$ if $\gamma_j = f_j(u_j)$, otherwise choose $\hat{u}_j\le u_j^*$, s.t.~$f_j(\hat{u}_j) = \gamma_j$. By construction of \eqref{FarjounSeibold:eq:condition_inequality_ingoing} and \eqref{FarjounSeibold:eq:condition_inequality_outgoing}, any resulting shocks and rarefaction waves are guaranteed to move away from the node (i.e.~forward on outgoing and backward in ingoing edges).

\section{Particleclaw}
\label{FarjounSeibold:sec:particleclaw}

\subsection{Characteristic Particles and Similarity Interpolant}
\label{FarjounSeibold:subsec:characteristic_particles_interpolation}
On each network edge, we have a scalar one-dimensional hyperbolic conservation law
\begin{equation}
\label{FarjounSeibold:eq:conservation_law}
u_t+\prn{f(u)}_x = 0\;,\quad
u(x,0) = u_0(x)\;,
\end{equation}
where the flux function $f$ is assumed to be twice differentiable and concave ($f''<0$) on the range of function values (see \cite{FarjounSeibold:FarjounSeibold2009} for extensions of the approach to flux functions with inflection points). We consider a special subset of exact solutions of \eqref{FarjounSeibold:eq:conservation_law}, which can be represented by a finite number of characteristic particles, as follows. A particle is a computational node that carries a position $x^i(t)$, and a function value $u^i(t)$. Note that we shall consistently denote particle indices by superscripts, while subscripts are reserved for edge indices. In the following, for convenience, we omit the time-dependence in the notation. Consider a time-dependent set of $n$ particles $P = \{(x^1,u^1),\dots,(x^n,u^n)\}$, where $x^1\le\dots\le x^n$. On the interval $[x^1,x^n]$ that is spanned by the particles, we define the \emph{similarity interpolant} $U_P(x)$ piecewise between neighboring particles as a true similarity solution of \eqref{FarjounSeibold:eq:conservation_law}. If $u^i \neq u^{i+1}$, it is implicitly given (and uniquely defined, since $f$ is concave) by
\begin{equation}
\label{FarjounSeibold:eq:interpolation}
\frac{x-x^i}{x^{i+1}-x^i} = \frac{f'(U_P(x))-f'(u^i)}{f'(u^{i+1})-f'(u^i)}\;.
\end{equation}
If $u^i = u^{i+1}$, the interpolant is simply constant $U_P(x) = u^i$. As shown in \cite{FarjounSeibold:FarjounSeibold2009}, the interpolant $U_P$ is an analytical solution of the conservation law \eqref{FarjounSeibold:eq:conservation_law}, as each particle moves according to the characteristic equations
\begin{equation}
\label{FarjounSeibold:eq:characteristic_equations}
\dot x = f'(u)\;,\quad
\dot u = 0\;,
\end{equation}
i.e.~$(x^i(t),u^i(t)) = (x^i(0)+f'(u^i(0))t,u^i(0))$. The reason for this fact is that due to the particular form of $U_P$, each point $(x(t),u(t))$ on it does move according to the same characteristic equations \eqref{FarjounSeibold:eq:characteristic_equations}. Note that strictly speaking, $U_P(x,t)$ is a solution only in the weak sense, since the derivative $\pd{}{x}U_P(x,t)$ is discontinuous at the particles. However, since $U_P$ is continuous, the Rankine-Hugoniot shock conditions \cite{FarjounSeibold:Evans1998} are trivially satisfied.

The interpolant $U_P$ is called ``similarity interpolant'', since the solution between neighboring particles \eqref{FarjounSeibold:eq:interpolation} is either a rarefaction wave that comes from a discontinuity (if $f'(x^i)<f'(x^{i+1})$) or a compression wave that will become a shock (if $f'(x^i)>f'(x^{i+1})$). As a consequence, the approach can be interpreted as \emph{rarefaction tracking} \cite{FarjounSeibold:FarjounSeibold2010}. This expresses both its similarities and differences to front tracking approaches \cite{FarjounSeibold:HoldenHoldenHeghKrohn1988, FarjounSeibold:HoldenRisebro2002}, which approximate the true solution by a finite number of shocks.

Just as the true solution of \eqref{FarjounSeibold:eq:conservation_law} may cease to be continuous after some critical time, the similarity approximation exists only up to the time of the first collision of neighboring particles. For a pair of neighboring particles $(x^i,u^i)$ and $(x^{i+1},u^{i+1})$, the time of collision (i.e.~they have the same $x$-position) is given by
\begin{equation}
\label{FarjounSeibold:eq:collision_time}
T^i = -\frac{x^{i+1}-x^i}{f'(u^{i+1})-f'(u^i)}\;,
\end{equation}
given that $f'(x^i)>f'(x^{i+1})$. Consequently, for $n$ particles, the time of the first collision is $T^* = \min\prn{\{T^i \ | \ T^i>0\}\cup\infty}$. At that time, a shock occurs (at $x^i = x^{i+1}$, from $u^{i+1}$ to $u^i$), and the method of characteristics cannot be used to evolve the solution further in time.

\subsection{Representation of Shocks}
\label{FarjounSeibold:subsec:shocks}
The characteristic motion of particles, as described in Sect.~\ref{FarjounSeibold:subsec:characteristic_particles_interpolation}, can only be performed if no shocks are present in the numerical solution. One idea to overcome this limitation, thus admitting an evolution of solutions with shocks, is to merge particles upon collision. This approach was first presented in \cite{FarjounSeibold:FarjounSeiboldMeshfree2008}, analyzed and generalized in \cite{FarjounSeibold:FarjounSeibold2009, FarjounSeibold:FarjounSeibold2010}, and implemented in the software \emph{particleclaw} \cite{FarjounSeibold:Particleclaw}. The merging of two particles $(x^i,u^i)$ and $(x^{i+1},u^{i+1})$ with $x^i = x^{i+1}$ into a single new particle $(x^i,\bar{u}^i)$ is performed such that the total area under the similarity interpolant $U_P$ is exactly preserved. As shown in \cite{FarjounSeibold:FarjounSeibold2009}, the area under $U_P$ between two neighboring particles equals
\begin{equation}
\label{FarjounSeibold:eq:area}
\int_{x^i}^{x^{i+1}}U_P(x)\,\ud{x} = (x^{i+1}-x^i)\,a(u^i,u^{i+1})\;,
\end{equation}
where
\begin{equation}
\label{FarjounSeibold:eq:average}
a(u^i,u^{i+1}) = \frac{\brk{f'(u)u-f(u)}_{u^i}^{u^{i+1}}}{\brk{f'(u)}_{u^i}^{u^{i+1}}}
\end{equation}
is a nonlinear average function (note that $\brk{f(u)}_v^w =
f(w)-f(v)$). Equating the area under the interpolant before and after
the merge, in general, yields a nonlinear equation for $\bar{u}^i$,
which can be solved up to machine accuracy by a few Newton iteration
steps (the geometry of the problem generally yields a very good
initial guess). In the case we are solving here, the interpolant is linear which simplifies many calculations.

The merging of colliding particles replaces a discontinuity by a
continuous interpolant, and thus the numerical approximation can be
evolved further in time using the characteristic particle motion
\eqref{FarjounSeibold:eq:characteristic_equations}. Since the merging
of particles $i$ and $i+1$ modifies the interpolant in the interval
$[x^{i-1},x^{i+2}]$, this approach introduces a small error right
around shocks. We control the magnitude of this error by the following
additional step. Let a parameter $d$ be given on the edge that
provides an upper bound on the distance of  particles adjacent to a pair of particles that need to be merged. Consider two particles $i$ and $i+1$ that need to be merged, because $x^i = x^{i+1}$ and $f'(u^i)>f'(u^{i+1})$. If $x^i-x^{i-1}>d$, then a new particle is inserted (on the interpolant) at $x^i-d$. Moreover, if $x^{x+2}-x^{x+1}>d$, then a new particle is inserted at $x^{i+1}+d$. After the relabeling of the particle indices, the original particles $i$ and $i+1$ are merged, in the way described above.

Note that there is an alternative to particle merges: the generalization to shock particles, as introduced in \cite{FarjounSeibold:FarjounSeiboldMeshfree2011}. The advantage of this version of the particle method is that shocks are represented exactly. The price to pay is a particle evolution that is more complicated than \eqref{FarjounSeibold:eq:characteristic_equations}. In this paper, for the application to macroscopic traffic models on networks, we consider the approach that uses characteristic particles. While it is certainly possible to formulate the ideas presented below with shock particles, the simplicity of the characteristic particle motion admits an easier presentation and analysis of the methodology. Moreover, the approach presented here introduces an intrinsic error at network vertices, and thus the exact nature of shock particles is less of an advantage.

\section{Generalizing Particleclaw to Network Flows}
\label{FarjounSeibold:sec:particleclaw_networks}
In this section we demonstrate how the approach \emph{particleclaw} can be generalized to nonlinear traffic models on highway networks. Our goal is to obtain an overall approach that is highly modular, in the sense that the solution on each edge can be evolved independently of the other edges, and edges communicate only during a \emph{synchronization step}. The specific methodology is as follows. Having sampled the initial conditions on each edge and synchronized them (as described below), we pick a time step $\Delta t$. During each time increment $t\in [t_n,t_n+\Delta t]$, the solution on each edge is evolved using the simple particle method described in Sect.~\ref{FarjounSeibold:sec:particleclaw}, with a special treatment at the first and the last particle (see Sect.~\ref{FarjounSeibold:subsec:virtual_domain_area_credit}). At the end of the time step, at each network node, the adjacent edges are synchronized with each other: first, the numerical solution (which may have partially moved away from the edge) is interpolated/extrapolated back onto the edge; second, the generalized Riemann problems (described in Sect.~\ref{FarjounSeibold:subsec:traffic_networks}) are invoked; lastly, area is suitably re-distributed. All these operations are performed such that the approach is exactly conservative, i.e.~no area is lost or created.

The independent evolution of the solution on the edges renders the approach extremely adapt towards parallelization: each edge can be stored in its own share of memory, communication between the different edges need only occur during the synchronization, and very little information must be transferred (see Sect.~\ref{FarjounSeibold:subsec:synchronization_step}). This methodology is possible due to the finite speed of information propagation of the hyperbolic conservation law \eqref{FarjounSeibold:eq:lighthill_whitham} on each edge. There is a maximum synchronization time step $\Delta t_\text{max}$, such that for all $\Delta t\le\Delta t_\text{max}$, information does not propagate further than half the length of each edge. As a consequence, in the synchronization step, all nodes can be treated independently of each other. Note that the maximum admissible time step $\Delta t_\text{max}$ between synchronization events is on the order of the smallest edge length divided by the fastest characteristic velocity. This is significantly larger than the maximum admissible time step in many traditional numerical approaches, which is on the order of the grid resolution divided by the fastest characteristic velocity. Hence with the presented particle approach, the relatively costly generalized Riemann problems need to be called much less frequently.

It should be pointed out that even though the method is exactly conservative, it is not exact. Due to the uncoupled evolution of the edges, a certain amount of information is lost, resulting in approximation errors that increase with the size of $\Delta t$. Thus, there could be accuracy constraints on $\Delta t$ that are more stringent than the stability constraints. Below, we outline the required additions to \emph{particleclaw} (Sect.~\ref{FarjounSeibold:subsec:virtual_domain_area_credit}), describe the synchronization step (Sect.~\ref{FarjounSeibold:subsec:synchronization_step}), and show how exact area balance is achieved (Sect.~\ref{FarjounSeibold:subsec:extrapolation}).

\subsection{Virtual Domain, Excess Area, Virtual Area, and Area Credit}
\label{FarjounSeibold:subsec:virtual_domain_area_credit}
One key idea of the particle approach is that the numerical solution can be advanced on each edge, without considering the coupling with other edges at the network nodes. Clearly, this evolution incurs an error near the two edge boundaries, since the coupling information is not used. However, in Sect.~\ref{FarjounSeibold:subsec:synchronization_step} we design the coupling step such that the arising gaps in information near the edge boundaries can be filled (with satisfactory accuracy) during the coupling step.

The coupling-free evolution on an edge $x\in [0,L]$ works very similarly to the basic \emph{particleclaw} methodology described in Sect.~\ref{FarjounSeibold:sec:particleclaw}, with a special twist at the first and the last particle, as described below. When no particles share the same $x$-position, then all particles are simply moved according to the method of characteristics \eqref{FarjounSeibold:eq:characteristic_equations}. Since this applies also to the first and the last particle, the domain of representation of the solution $[x^1(t),x^N(t)]$ does in general not match with the computational domain $[0,L]$. In the case that particles extend beyond the edge into the \emph{virtual domain}, i.e.~$x^1<0$ or $x^N>L$, the similarity interpolant \eqref{FarjounSeibold:eq:interpolation} defines the numerical solution in particular up to $x=0$ or $x=L$, respectively. In addition, it defines an \emph{excess area} $\int_{x^1}^0 U_P(x)\,\ud{x}$ or $\int_L^{x^N} U_P(x)\,\ud{x}$, respectively. In the case that particles do not reach the edge boundary, i.e.~$x^1>0$ or $x^N<L$, we can define an extrapolation on $[0,x^1]$ or $[x^N,L]$, respectively, if we are given a value for the \emph{virtual area} $\int_0^{x^1} u\,\ud{x}$ or $\int_{x^N}^L u\,\ud{x}$, respectively. How to construct this extrapolation is described in Sect.~\ref{FarjounSeibold:subsec:extrapolation}.

The merging of colliding particles works as described in Sect.~\ref{FarjounSeibold:subsec:shocks}, with one exception. If for the merging of particles $i$ and $i+1$, there is no particle $i-1$, then a new particle is added at $x^i-d$ with the same value as $u^i$, and a ``left-sided'' \emph{area credit} $I_L = d\,u^i$ is recorded. Similarly, if there is no particle $i+2$, a new particle is added at $x^{i+1}+d$ with the same value as $u^{i+1}$, and a ``right-sided'' \emph{area credit} $I_R = d\,u^{i+1}$ is recorded. After these insertions, the merge can be performed.

\subsection{Synchronization Step}
\label{FarjounSeibold:subsec:synchronization_step}
The synchronization of area between the edges happens by moving area
from edges in which information is flowing into the node to edges from which information is flowing out of the node. As introduced in Sect.~\ref{FarjounSeibold:subsec:traffic_networks}, we denote the edges on which vehicles enter the node ``ingoing''. On the ingoing edge $i\in\{1,\dots,n\}$, the position of the node is at $x = L_i$. Conversely, the edges on which vehicles exit the node are called ``outgoing''. On the outgoing edge $j\in\{n+1,\dots,n+m\}$, the position of the node is at $x = 0$. Since the LWR model \eqref{FarjounSeibold:eq:lighthill_whitham} admits characteristic velocities of either sign, information can propagate either into or away from the node, both for ingoing and outgoing edges. Therefore, we introduce further terminology: edges for which information is going into the node are called ``influencing'', edges for which information is going away from the node are denoted ``affected'', and edges with zero characteristic velocity are called ``neutral''. When solving a generalized Riemann problem at a node, the procedure described in the last paragraph of Sect.~\ref{FarjounSeibold:subsec:traffic_networks} implies that all edges for which the flux changes (from $f_i(u_i)$ to $\gamma_i$) become automatically affected or neutral. Influencing edges arise if $\gamma_i = f_i(u_i)$ and in addition: $u_i<u_i^*$ for ingoing edges and $u_i>u_i^*$ for outgoing edges.

Let us now focus on one node in the network, with $n$ ingoing and $m$ outgoing edges. In the following, subscripts denote the edge index, and superscripts denote the particle index. With this notation, the last particle on the ingoing edge $i\in\{1,\dots,n\}$ is $(x_i^{N_i},u_i^{N_i})$, and the first particle on the outgoing edge $j\in\{n+1,\dots,n+m\}$ is $(x_j^1,u_j^1)$. We design our approach such that at the end of the synchronization step (and thus at the beginning of the evolution on each edge) the first and last particle on each edge arise as solutions of the generalized Riemann problems described in Sect.~\ref{FarjounSeibold:subsec:traffic_networks}. Consequently, due to \eqref{FarjounSeibold:eq:driver_destinations}, we have the $m$ conditions
\begin{equation}
\label{FarjounSeibold:eq:flux_relations_destinations}
\sum_{i=1}^n a_{i,j} f_i(u_i^{N_i}) = f_j(u_j^1) \quad\forall\,j\in\{n+1,\dots,n+m\}\;,
\end{equation}
given by the desired destinations matrix $A$, and---if necessary---further $n-1$ conditions
\begin{equation}
\label{FarjounSeibold:eq:flux_relations_merging}
f_i(u_i^{N_i}) = \beta c_i \quad\forall\,i\in\{1,\dots,n\}
\end{equation}
for some $\beta\in\mathbb{R}$, given by the merging vector $\vec{c}\in\mathbb{R}^n$.

We now derive, first for ingoing edges, an evolution equation for the excess area or virtual area (see Sect.~\ref{FarjounSeibold:subsec:virtual_domain_area_credit}) that is between the end of the edge $L_i$ and the last particle $x_i^{N_i}$. Letting $u$ denote a true solution of the hyperbolic conservation law \eqref{FarjounSeibold:eq:lighthill_whitham}, this area $I_i = \int_{L_i}^{x_i^{N_i}} u(x,t)\,\ud{x}$ has the rate of change
\begin{equation}
\label{FarjounSeibold:eq:evolution_virtual_area_ingoing}
\frac{d}{dt} I_i = f_i(u_i(L_i,t))-f_i(u_i^{N_i})+u_i^{N_i}f_i'(u_i^{N_i})\;.
\end{equation}
Here, the first two terms come from the PDE \eqref{FarjounSeibold:eq:lighthill_whitham}, and the third term stems from the fact that $x_i^{N_i}$ moves with velocity $f_i'(u_i^{N_i})$. Analogously, we compute the evolution of the excess/virtual area $I_j = \int_{0}^{x_1^0} u(x,t)\,\ud{x}$ on an outgoing edge as
\begin{equation}
\label{FarjounSeibold:eq:evolution_virtual_area_outgoing}
\frac{d}{dt} I_j = f_j(u_j(0,t))-f_j(u_j^1)+u_j^1 f_j'(u_j^1)\;.
\end{equation}

At the beginning of the coupling-free evolution, the first and last particle are on the edge boundaries, and thus all excess/virtual areas are zero. At the end of the coupling-free evolution, at the considered node, we know the excess area for all influencing edges. Moreover, for all neutral edges the excess area is zero. The idea is now to use the evolution equations \eqref{FarjounSeibold:eq:evolution_virtual_area_ingoing} and \eqref{FarjounSeibold:eq:evolution_virtual_area_outgoing}, together with the relations \eqref{FarjounSeibold:eq:flux_relations_destinations} and \eqref{FarjounSeibold:eq:flux_relations_merging}, to determine the unknown virtual areas at all affected edges. Any area credit that has been taken due to merges at the first or last particle one an edge (see Sect.~\ref{FarjounSeibold:subsec:virtual_domain_area_credit}) is simply subtracted from this balance. Below, we present this idea for three fundamental cases: a bottleneck, a bifurcation, and a confluence. Many practically relevant road networks can be constructed from these three types of nodes.

\subsubsection{Bottleneck}
A bottleneck is a sudden change of the road conditions, i.e.~$n=1$ and $m=1$. Here, the feasibility domain \eqref{FarjounSeibold:eq:feasibility_domain} is $\Omega = [0,d_1]\cap [0,s_2]$, where $d_1 = f_1(\min\{u_1,u_1^*\})$ and $s_2 = f_2(\max\{u_2,u_2^*\})$. Thus, the solution of the optimization problem \eqref{FarjounSeibold:eq:node_optimization} is $\gamma_1 = \min\{d_1,s_2\}$, which allows for two cases. If $\gamma_1 = d_1$, then all vehicles can pass through the node, and consequently edge 1 is influencing (or neutral), and edge 2 is affected (or neutral). Otherwise, if $\gamma_1<d_1$, then not all vehicles can pass, and a jam is triggered. Consequently, edge 1 is affected, and edge 2 is influencing or neutral.

\begin{figure}[t]
\includegraphics[width=\textwidth]{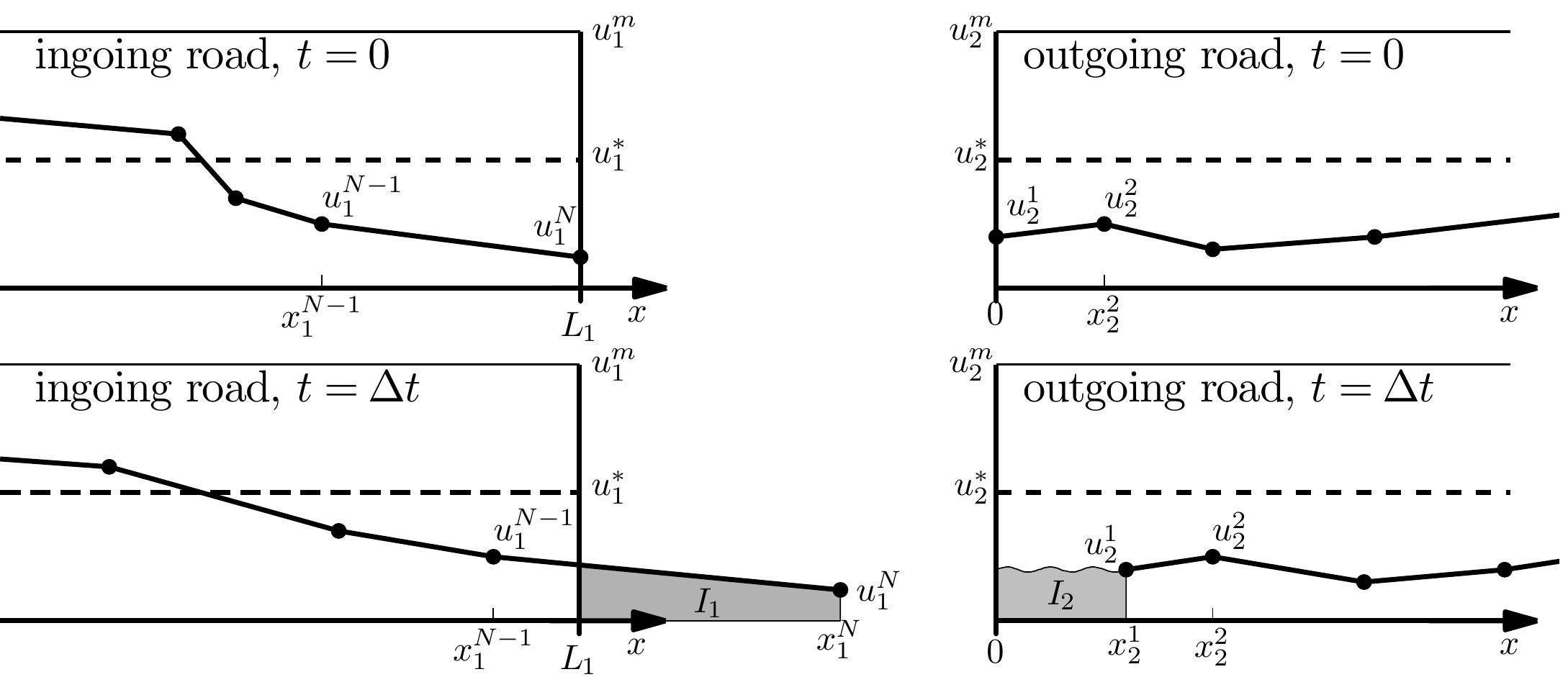}
\caption{Area tracking for a bottleneck. The excess area $I_1$ allows the computation of the virtual area $I_2$.}
\label{FarjounSeibold:fig:bottleneck_area}
\end{figure}

In either case, since the influx equals the outflux at all times, we have that $f_1(u_1(L_1,t)) = f_2(u_2(0,t))$ and $f_1(u_1^{N_1}) = f_2(u_2^1)$, and thus the evolution equations \eqref{FarjounSeibold:eq:evolution_virtual_area_ingoing} and \eqref{FarjounSeibold:eq:evolution_virtual_area_outgoing} imply that
\begin{equation*}
\frac{d}{dt} (I_2-I_1) = u_2^1 f_2'(u_2^1)-u_1^{N_1}f_1'(u_1^{N_1}) = C\;,
\end{equation*}
which is a known time-independent quantity. Together with $I_1(t) = I_2(t) = 0$, we obtain a relation $I_2(t+\Delta t)-I_1(t+\Delta t) = C\Delta t$ that allows to obtain $I_2$ from $I_1$ (in the case $\gamma_1 = d_1$), or $I_1$ from $I_2$ (in the case $\gamma_1<d_1$). Figure~\ref{FarjounSeibold:fig:bottleneck_area} shows an example of the movement of particles into the virtual domain, and the evolution of excess area ($I_1$) and virtual area ($I_2$).

\subsubsection{Bifurcation}
In a bifurcation, one road splits into two, i.e.~$n=1$ and $m=2$. The relative numbers of vehicles that exit the node on edge 2 or edge 3 are given by $a_{1,2}$ and $a_{1,3} = 1-a_{1,2}$, respectively. Here, the feasibility domain \eqref{FarjounSeibold:eq:feasibility_domain} is $\Omega = \{\gamma_1\in [0,d_1]\ |\ a_{1,j}\gamma_1\in [0,s_j]\;\forall\,j\in\{2,3\}$, where $d_1 = f_1(\min\{u_1,u_1^*\})$ and $s_j = f_j(\max\{u_j,u_j^*\})\;\forall\,j\in\{2,3\}$. The solution of the optimization problem \eqref{FarjounSeibold:eq:node_optimization} is given by $\gamma_1 = \min\{d_1,\frac{s_2}{a_{1,2}},\frac{s_3}{a_{1,3}}\}$. This allows two possibilities. If $\gamma_1 = d_1$, then all vehicles can pass through the node, and edge 1 is influencing (or neutral), and edges 2 and 3 are affected (or neutral). Otherwise, if $\gamma_1<d_1$, then not all vehicles can pass, and edge 1 becomes affected. Of the outgoing edges, the one that causes the congestion is influencing or neutral, while the other outgoing edge is affected or neutral.

In either situation, we can use the evolution equations \eqref{FarjounSeibold:eq:evolution_virtual_area_ingoing} and \eqref{FarjounSeibold:eq:evolution_virtual_area_outgoing} to determine the virtual areas, since \eqref{FarjounSeibold:eq:flux_relations_destinations} provides two conditions that relate $I_1$, $I_2$, and $I_3$ with each other. As an example, consider the case of edge 2 being influencing, and edges 1 and 3 affected. Here, $I_2$ is known (by integrating the interpolant \eqref{FarjounSeibold:eq:interpolation}). From \eqref{FarjounSeibold:eq:evolution_virtual_area_ingoing} and \eqref{FarjounSeibold:eq:evolution_virtual_area_outgoing} we first compute $C = \frac{d}{dt} (I_1-\tfrac{1}{a_{1,2}}I_2) = u_1^{N_1}f_1'(u_1^{N_1})-\tfrac{1}{a_{1,2}}u_2^1 f_2'(u_2^1)$, and then obtain $I_1(t+\Delta t) = \tfrac{1}{a_{1,2}}I_2(t+\Delta t)+C\Delta t$. Knowing $I_1$, we can then determine $I_3$, by considering $\frac{d}{dt} (I_3-a_{1,3}I_1)$.

\subsubsection{Confluence}
Here, two roads converge into one, i.e.~$n=2$ and $m=1$. Flux conservation simply states that $\gamma_1+\gamma_2 = \gamma_3$. The feasibility domain \eqref{FarjounSeibold:eq:feasibility_domain} is now two-dimensional: $\Omega = \{(\gamma_1,\gamma_2)\in [0,d_1]\times [0,d_2] \ | \ \gamma_1+\gamma_2\in [0,s_3]\}$, where $d_i = f_i(\min\{u_i,u_i^*\})\;\forall\,i\in\{1,2\}$ and $s_3 = f_3(\max\{u_3,u_3^*\})$. For the optimization problem \eqref{FarjounSeibold:eq:node_optimization}, we distinguish two cases. If $d_1+d_2\le s_3$, then all vehicles can pass through the node, and the ingoing edges are influencing (or neutral), and the outgoing edge is affected (or neutral). In contrast, if $d_1+d_2>s_3$, then a jam is triggered. How the backwards going information distributes among the two ingoing edges depends now on the merging vector $(c_1,c_2)$. The following outcomes are possible: (i) all vehicles from edge 1 get through, i.e.~$\gamma_1 = d_1$ and a jam occurs on edge 2, i.e.~$\gamma_2<d_2$. In this case, edge 1 is influencing, edge 2 is affected, and edge 3 is influencing or neutral; (ii) a jam occurs on both ingoing edges, i.e.~$\gamma_1<d_1$ and $\gamma_2<d_2$. Here, both ingoing edges are affected, and the outgoing edge is influencing or neutral; (iii) is the same as (i), but with edges 1 and 2 reversed.

Again, the virtual areas at the affected edges can be obtained systematically from the known virtual areas at the influencing and neutral edges. In the case that all vehicles pass through the node, we have $I_3(t+\Delta t) = I_1(t+\Delta t)+I_2(t+\Delta t)+\prn{u_3^1 f_3'(u_3^1)-u_1^{N_1}f_1'(u_1^{N_1})-u_2^{N_2}f_2'(u_2^{N_2})}\Delta t$. In turn, if jamming occurs, e.g.~on both ingoing edges, we use the knowledge of the merging ratios $\frac{\gamma_1}{\gamma_2} = \frac{c_1}{c_2}$ to determine $I_1$ and $I_2$ from $I_3$. The methodology is analogous to the cases presented above.

\subsection{Representation of Virtual Area by Particles}
\label{FarjounSeibold:subsec:extrapolation}
In Sect.~\ref{FarjounSeibold:subsec:synchronization_step}, we have described how the virtual area at affected edges can be determined. The presented methodology tracks area exactly, and thus the numerical approach is exactly conservative. In order to finalize the synchronization step, we must account for the imbalance of area by changing parts of the solution, thereby creating or removing area.

Throughout the design of \emph{particleclaw}, the two main guiding principles are: (i) conserve $\int u\,\ud{x}$ exactly and as locally as possible; and (ii) ensure that the approach is TVD, i.e.~no spurious overshoots or undershoots are generated. We conclude the synchronization step by modifying the numerical solution, while following these two principles. The conservation principle means that the excess area and the virtual area must be converted into actual area under particles on the edges near where that area was recorded. The locality requirement is respected by a finite domain of influence: the solution is never modified further than $\Delta t$ times the maximal particle velocity from the node. The TVD condition means that the modified solution does not possess values above or below the range of values that is presented in the nearby particles and the solution $\hat u$ that comes from the generalized Riemann problem~\eqref{FarjounSeibold:eq:node_optimization}. The implementation outlined below violates the TVD condition only under exceptional circumstances. In fact, in the numerical experiments presented in Sect.~\ref{FarjounSeibold:sec:numerical_results}, this non-TVD ``last resort'' solution turned out to never be required.

The precise methodology is as follows. First, the numerical solution is truncated or extrapolated to the edges, by adding a particle either on the interpolant, or with the same value as the closest particle, respectively. All excess areas (i.e.~the parts that extend beyond the domain of each edge) can now be removed, since their contribution has been accounted for in the step presented in Sect.~\ref{FarjounSeibold:subsec:synchronization_step}. Now, the generalized Riemann solver is called with the new values on the node position. The resulting $\hat u$-values are added to the solution as new extremal particles (unless the value has not changed, $\hat u = u$, in which case no particle is added) with the same $x$-value as the extremal particle. Being able to have two particles with the same $x$-value is a peculiarity of the particle method, and it implies that the particle insertion does not change the area under the solution.

\begin{figure}[t]
\includegraphics[width=\textwidth]{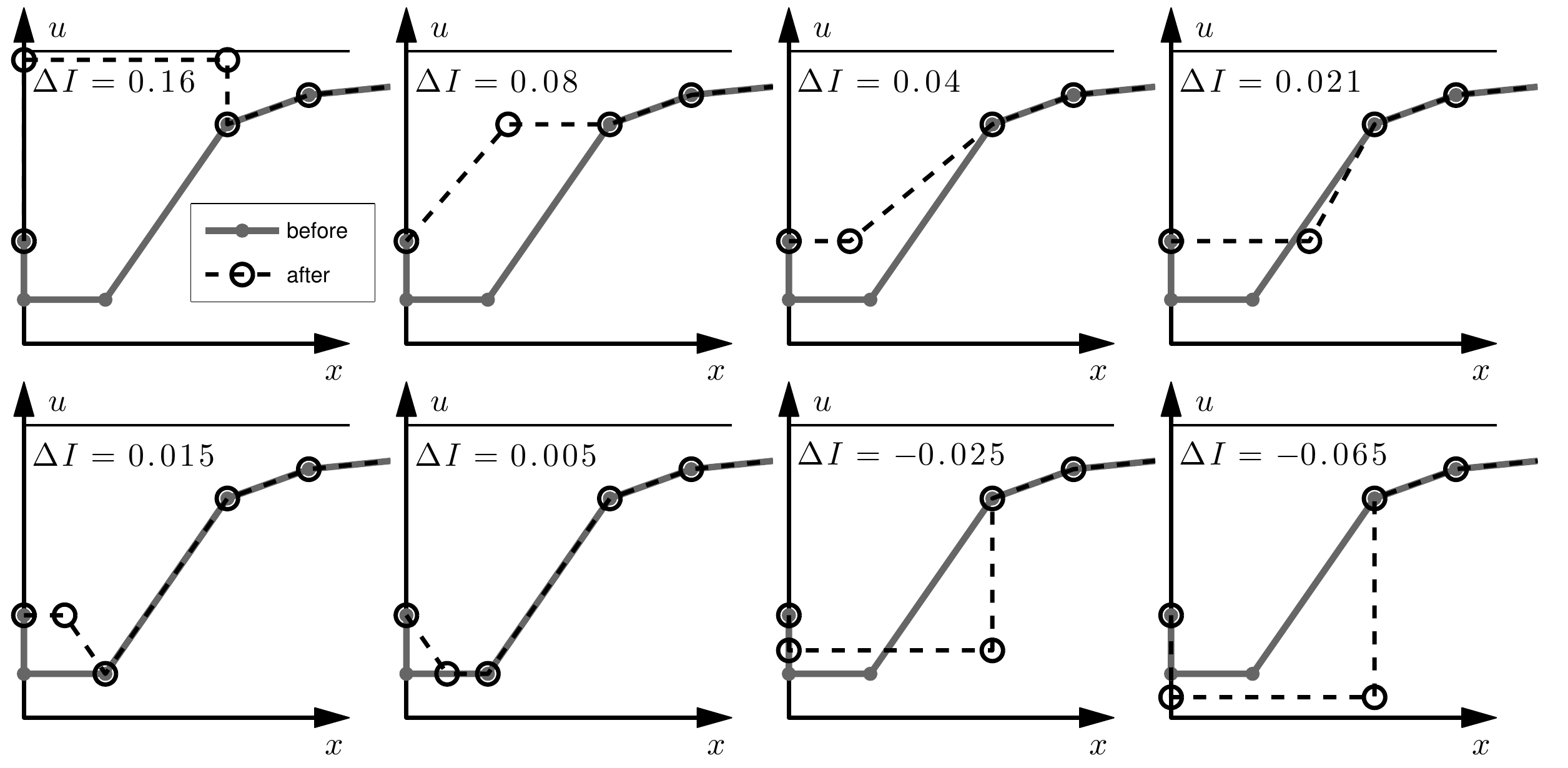}
\caption{Reconstructed solution for various $\Delta I$. The thick gray line is the original solution, and the dashed black line the ``correct area'' solution that replaces it, whose area is greater by $\Delta I$.}
\label{FarjounSeibold:fig:area:reconst}
\end{figure}

With this, on each edge we have a numerical solution, whose extremal states are solutions of the generalized Riemann problems. What is left is the accounting for the correct balance of area. We do so by modifying the numerical solution, however, without changing the previously found new extremal states. This guarantees that at the end of the synchronization step, the fluxes at the all the nodes are solutions of the generalized Riemann problem. In the preceding subsections we have derived the area $I$ that is supposed to be under the solution between the network node and the extremal particle. In general, the area under the numerical solution constructed above by simple interpolation/extrapolation has a small disparity with $I$. The difference, $\Delta I$, needs to be added to the affected edge, by modifying the numerical solution. Without loss of generality let us examine an \emph{outgoing} affected edge. We look for a new solution to replace the current solution near the end of the edge. This new solution must have an area difference $\Delta I$ from the area of the current solution, and the first particle must remain unchanged in order to satisfy the generalized Riemann problem. We follow the steps below (starting with $k=2$) to represent the required change in area using particles (illustrated in Fig.~\ref{FarjounSeibold:fig:area:reconst}):
\begin{itemize}
\item[(1)] Attempt a solution that comprises a constant part and a
  linear part so that it connects particles $1$ and $k$ in a
  continuous fashion. If such a solution can provide $\Delta I$, accept it.
\item[(2)] Attempt a solution with constant value $u$ in the
  domain $[x_1,x_k]$, where $u\in [\min_{i\le k}(u_i),\allowbreak\max_{i\le
    k}(u_i)]$. If such a solution can provide $\Delta I$, accept it. If this solution is not accepted, and $x_k<\Delta t\max_u(f'(u))$, increase $k$ by one and go back to (1), otherwise, go to (3).
\item[(3)] Use a solution that has a constant value $u$ in the domain $[x_1,x_k]$ with the value of $u$ needed to match $\Delta I$. 
\end{itemize}
In the algorithm above, step (3) is a fail-safe mechanism that will always result in a solution, but the solution may be non-TVD or even unphysical.

Figure~\ref{FarjounSeibold:fig:area:reconst} shows eight exemplary cases using the same initial particle configuration but different $\Delta I$.
In these examples, $k$ could not increase beyond 4: 
$\Delta I=0.16$ resulted in a solution at step (3) with $k=4$; $\Delta I=0.08,\,0.04,\,0.021$, with step (1) and $k=4$; $\Delta I=0.015,\, 0.005$ with step (1) and $k=3$; $\Delta I=-0.025$  with step (2) and $k=4$; and $\Delta I=-0.065$ with step (3) and $k=4$.

Performing this area reconstruction on all affected edges guarantees that area is preserved exactly; it may increase or decrease the number of particles.

\section{Numerical Results}
\label{FarjounSeibold:sec:numerical_results}

\subsection{Bottleneck Test Case}
\begin{figure}[t]
\begin{minipage}[b]{.52\textwidth}
\includegraphics[width=1.03\textwidth]{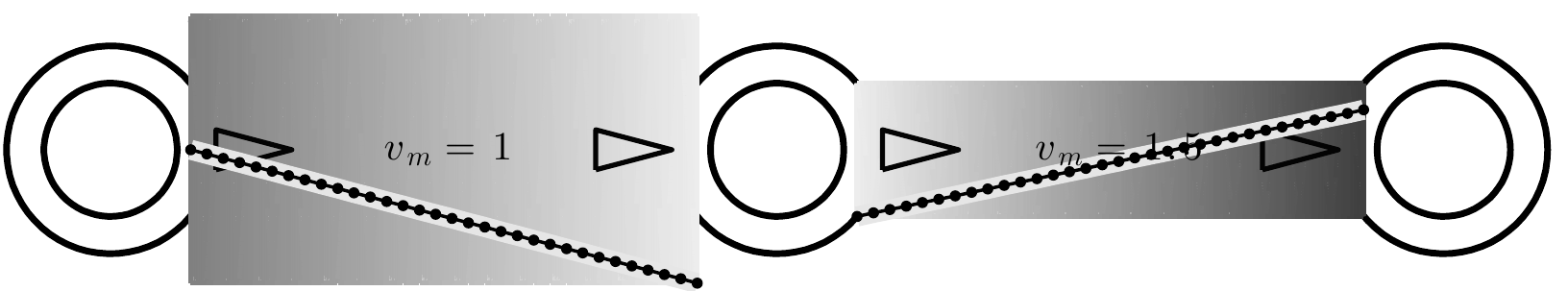} \\[3em]
\includegraphics[width=1.03\textwidth]{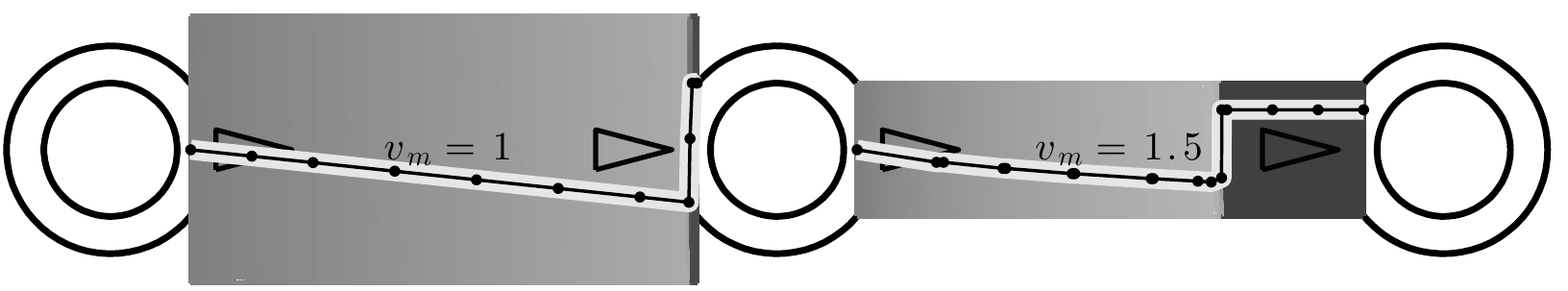} \\[3em]
\includegraphics[width=1.03\textwidth]{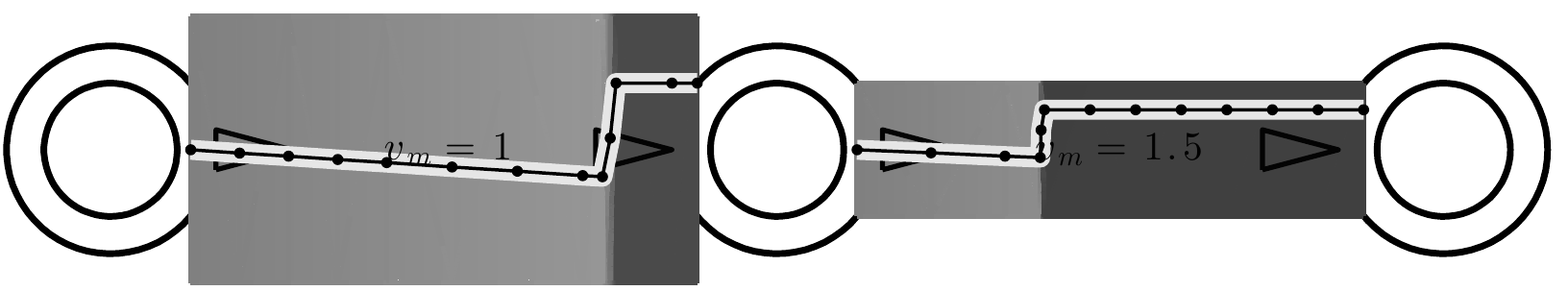} \\[-.2em]
\caption{Bottleneck test case. Solution at time $t=0$ (top), $t=1.5$ (middle), and $t=3$ (bottom), computed with \emph{particleclaw} with $h = 8\cdot 10^{-2}$ and $d = 2\cdot 10^{-2}$.}
\label{FarjounSeibold:fig:bottleneck_solution}
\end{minipage}
\hfill
\begin{minipage}[b]{.45\textwidth}
\hfill\includegraphics[width=\textwidth]{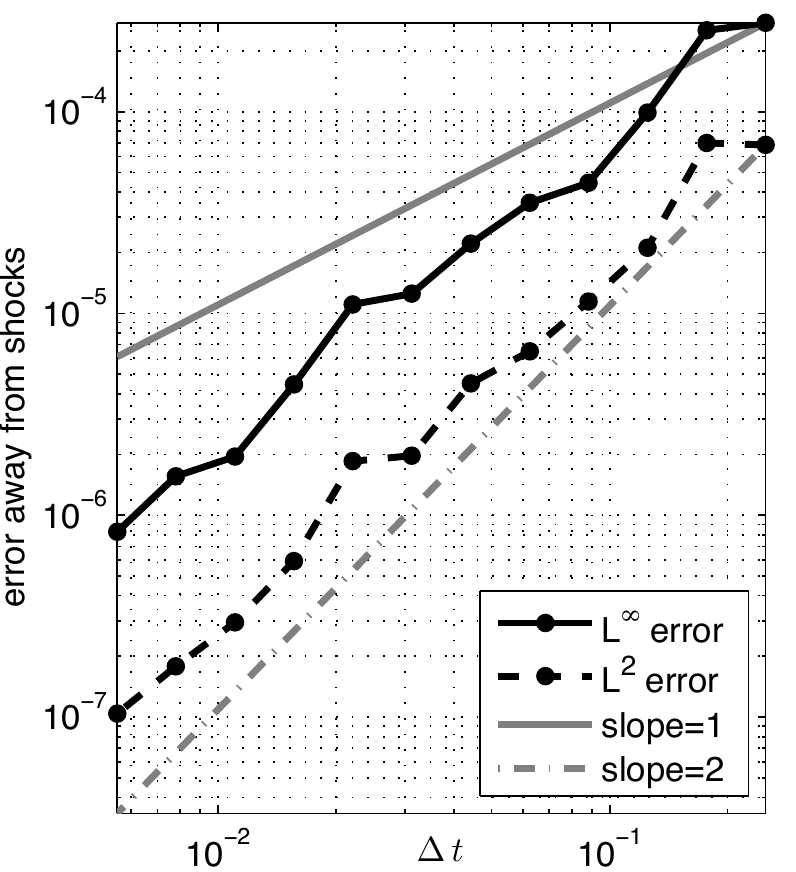} \\[-1.6em]
\caption{Bottleneck test case. The error as a function of the synchronization time step $\Delta t$ exhibits second order convergence behavior.}
\label{FarjounSeibold:fig:bottleneck_error}
\end{minipage}
\end{figure}

\begin{figure}[t]
\includegraphics[width=.9\textwidth]{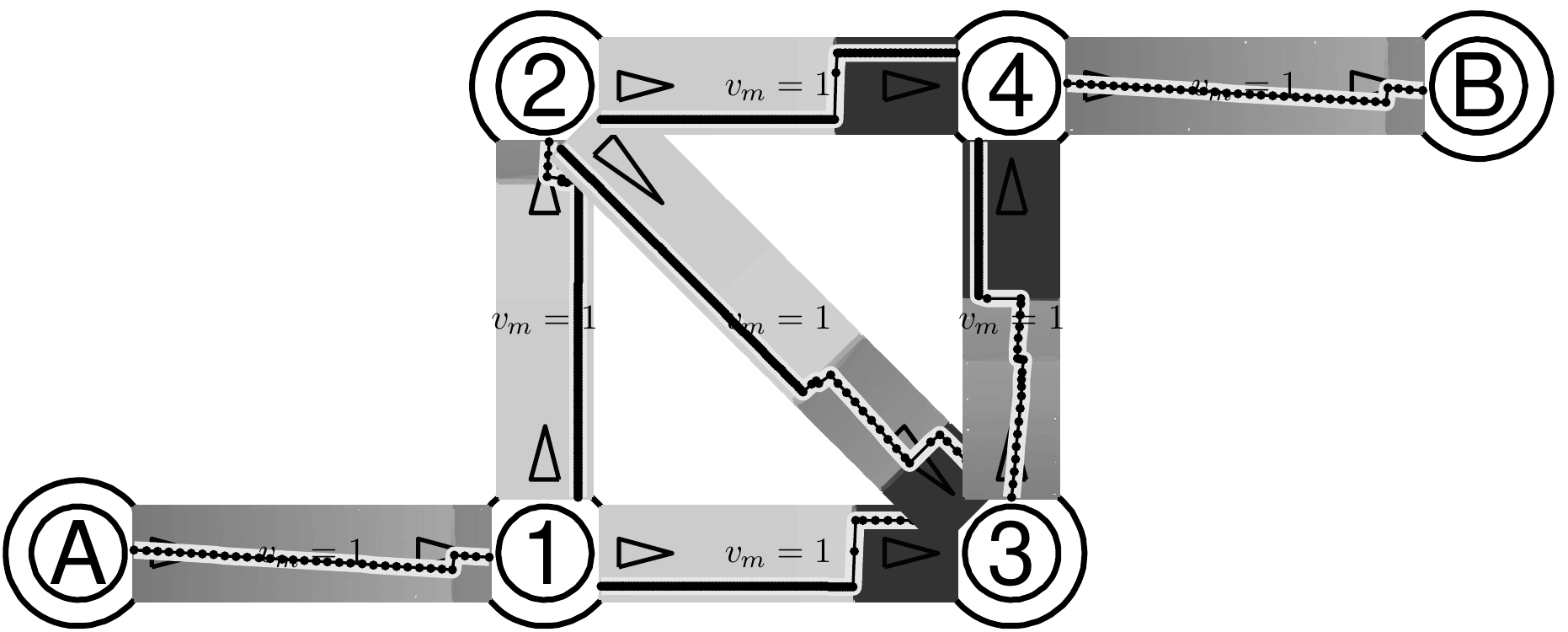} \\[1em]
\includegraphics[width=.9\textwidth]{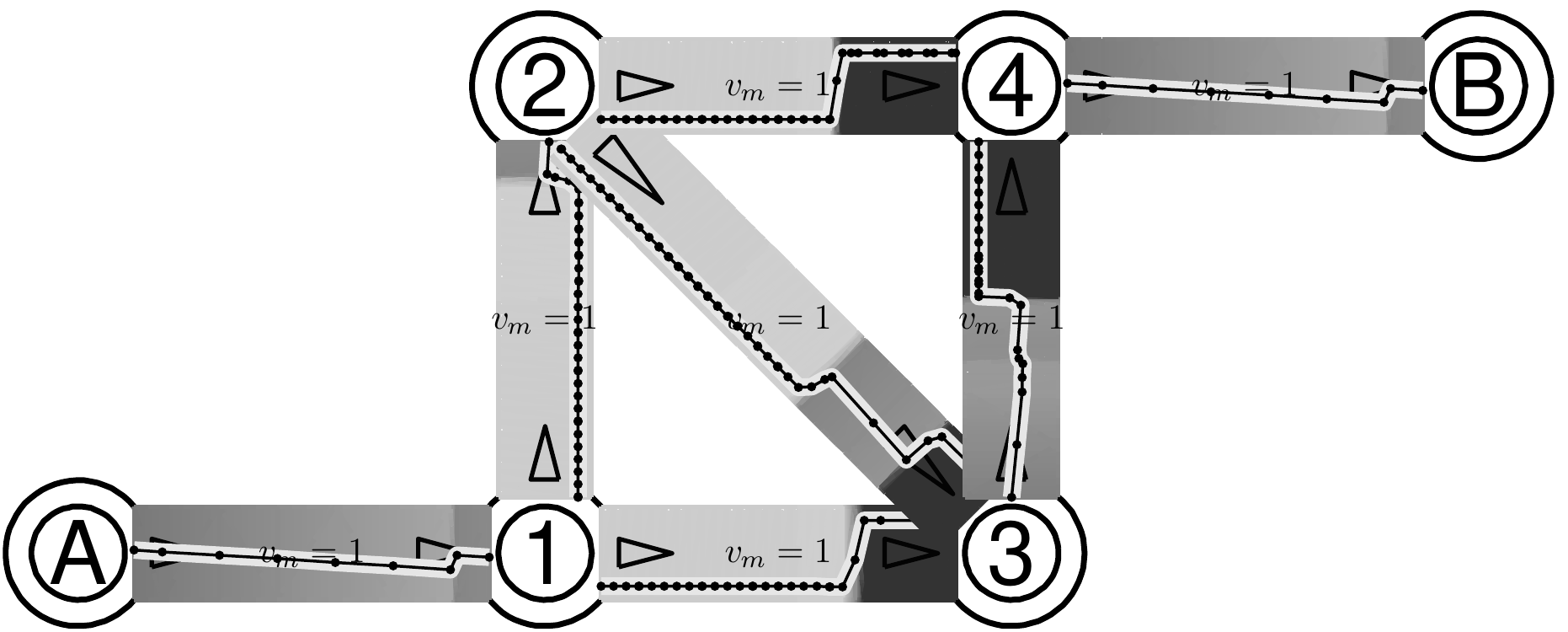} \\[1em]
\includegraphics[width=.9\textwidth]{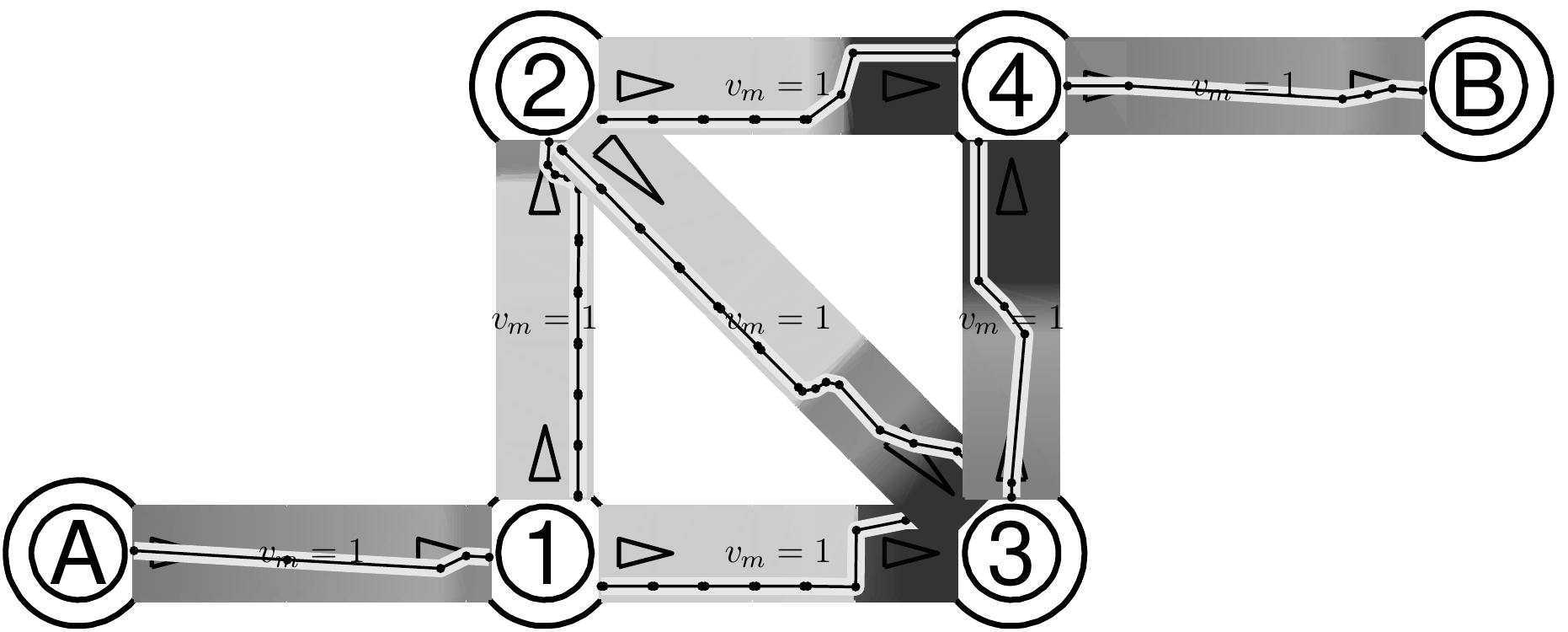}
\caption{Diamond network test case with various resolutions. The solutions have parameters $h = 2\cdot 10^{-2} s $, $d = 5\cdot 10^{-3} s $, and $\Delta t = 1\cdot 10^{-2} s$, where $s=1$ at the top, $s=5$ in the middle, and $s=20$ at the bottom.}
\label{FarjounSeibold:fig:diamond_network}
\end{figure}

As a first test case, we consider the simplest non-trivial network: a bottleneck that consists of two edges, each of length $L_1 = L_2 = 1$ that are joined linearly, as shown in Fig.~\ref{FarjounSeibold:fig:bottleneck_solution}. In dimensionless variables, the maximum traffic densities on the two roads are $u_1^\text{m} = 2$ and $u_2^\text{m} = 1$, and the maximum vehicle velocities are $v_1^\text{m} = 1$ and $v_2^\text{m} = 1.5$. This example can be interpreted as a model for the situation of two lanes merging into one lane, and at the same position the speed limit increasing by $50\%$. The initial conditions are $u_1(x) = 1-x$ and $u_2(x) = 0.8\,x$, and the boundary conditions of $u_1(0,t) = 1$ and $u_2(1,t) = 0.8$ are prescribed whenever information is entering the edge through these boundaries. The final time is $t_\text{final} = 3$. In Fig.~\ref{FarjounSeibold:fig:bottleneck_solution} (and also in Fig.~\ref{FarjounSeibold:fig:diamond_network}), the solution is shown in two ways: the shade of gray becomes darker with higher values of $u$, and the plot of $u$ vs.~$x$ is overlaid with dots representing the particles of the method. Each edge is annotated with the maximum velocity $v_m$ and two arrows showing the direction of vehicles. The thickness of each segment is proportional to its maximum density $u_m$.

It is shown in \cite{FarjounSeibold:FarjounSeibold2009} that \emph{particleclaw} itself is second-order accurate with respect to the maximum distance of particles (away from shocks or when using shock sharpening). What we investigate here is the error as a function of the size of the synchronization time step $\Delta t$, while having a large number of particles, such that the error due to the spacial approximation is negligible. Specifically, we provide the following two parameters for \emph{particleclaw}: an initial/desired distance of particles of $h = 8\cdot 10^{-5}$ is given, and the distance of particles that are inserted near a shock is $d = 2\cdot 10^{-5}$. Both parameters result in a very fine resolution, using more than ten thousand particles on each edge. We choose the time step from $\Delta t = 2^{-k/2}$, where $k\in\{4,\dots,15\}$, and compute a reference solution using $k=20$. Since here we do not wish to measure the error at a shock, we evaluate the error in the segment $x\in [0,0.3]$ on edge 2, in which the solution is smooth, but non-linear. Figure~\ref{FarjounSeibold:fig:bottleneck_error} shows the $L^\infty([0,\,0.3])$ and $L^2([0,\,0.3])$ errors of the approximate solution as a function of the synchronization time step, $\Delta t$. One can see a second-order dependence, i.e.~the error scales like $O((\Delta t)^2)$.

\subsection{Simulation of a Diamond Network}
An important point of \emph{particleclaw} is that it tends to yield
high quality numerical solutions (including shocks) for very few
particles (see e.g.~\cite{FarjounSeibold:FarjounSeibold2009,
  FarjounSeibold:FarjounSeiboldMeshfree2011}). Here, we demonstrate this feature for
a so-called diamond network. This network of seven edges and six nodes
is shown in Fig.~\ref{FarjounSeibold:fig:diamond_network}. Vehicles enter in the
node marked ``A'' driving towards ``1''; then the traffic flow splits
into a N and a E direction. At ``2'', drivers again decide to go SE or
straight E. At ``3'' and ``4'', two confluences happen, and finally
vehicles exit the network at ``B''. We choose all roads identical
$u_1^\text{m} = 1$ and $v_1^\text{m} = 1$, and at all network nodes,
an even split of traffic flow occurs. The lengths of all the roads is
taken to be 1 (even if it does not seem so in the plot). The initial
conditions are also identical on all the roads:
$u(x,0)=0.4+0.4\cos(3\pi x)$. We evolve the solution until the time
$t_\text{final} = 2$. Figure~\ref{FarjounSeibold:fig:diamond_network} shows the
computational results when applying \emph{particleclaw}, with three
different types of resolutions (in space and time). The top solution
has an initial/desired distance of particles of $h = 2\cdot 10^{-2}$,
the distance of particles that are inserted near a shock is $d =
5\cdot 10^{-3}$, and the synchronization time step is $\Delta t =
1\cdot 10^{-2}$. In the middle solution, all of these numbers are
multiplied by a factor of 5, and in the bottom solution by a
factor of 20. That means that the top possesses 20 times the
resolution than the bottom, and it is therefore about 400 times as
costly to compute. One can observe that the higher resolution at the
top yields sharper features. However, the presence and position of
the shocks (and other features) is well captured by the resolution in
the middle, and even at the bottom most features are visible quite clearly.

\section{Conclusions and Outlook}
\label{FarjounSeibold:sec:conclusions_outlook}
The framework presented in this paper demonstrates that the characteristic particle method \emph{particleclaw} can be generalized to solving first order traffic models on road networks in a robust, accurate, and exactly conservative fashion. The presented methodology allows to apply the basic one-dimensional solver on each edge. The coupling of traffic states at network nodes is achieved by a special synchronization step that ensures exact conservation properties by applying a proper transfer of area between edges. A fundamental ingredient in this synchronization is the invoking of generalized Riemann solvers that have been developed for traffic networks. In this paper, a specific focus lies on 1-to-1, 1-to-2, and 2-to-1 highway junctions. However, the presented methodology applies to more general types of networks, as long as the evolution on each edge is described by a scalar hyperbolic conservation law, and the problem description allows to formulate coupling conditions at the network nodes.

As previously shown \cite{FarjounSeibold:FarjounSeibold2009}, \emph{particleclaw} itself can be made second-order accurate with respect to the spacing between particles. The numerical examples investigated here confirm that the presented method is also second order accurate with respect to the time step $\Delta t$ between synchronization events. Moreover, it is demonstrated that the presented approach yields good quality numerical solutions (with an accurate location of shocks), even when using only few particles on each edge.

In contrast to many traditional numerical methods, here the synchronization time step $\Delta t$ is not limited by the spacing between particles (i.e.~no CFL stability restriction). It is only limited by the length of the edges, insofar as between two synchronization events, information is not allowed to travel beyond a single edge. In general, this allows for much larger values for $\Delta t$ than traditional methods do. An exception is posed by networks that possess a few very short edges. These can be treated significantly more efficiently by a simple generalization of the presented approach, namely by using a different $\Delta t$ for each network node, depending on the shortest edge connected to it. Thus, edges that neighbor short edges would need to synchronize more often, but the rest of the network would not.

The accurate location of shocks with only few particles on each edge, and the possibility of invoking the coupling conditions at network nodes relatively rarely, allow for fast and memory-efficient implementations. This makes the approach attractive for the simulation of large road networks, as well as the optimization, design, and control of the traffic flow on such large networks. It is also conceivable that the presented methodology can be adapted to other kinds of nonlinear flows on networks, for instance continuum models of supply chains \cite{FarjounSeibold:ArmbrusterMarthalerRinghofer2004, FarjounSeibold:HertyKlarPiccoli2007}.

\section*{Acknowledgments}
Y. Farjoun was financed by the Spanish Ministry of Science and Innovation under grant FIS2008-04921-C02-01. The authors would like to acknowledge support by the National Science Foundation. B. Seibold was supported through grant DMS--1007899. In addition, B. Seibold wishes to acknowledge partial support by the National Science Foundation through grants DMS--0813648 and DMS--1115278. Y. Farjoun wishes to acknowledge partial support through grant DMS--0703937.

\bibliographystyle{siam}
\bibliography{SeiboldFarjounReferences}

\end{document}